\documentclass[11pt]{article}
    \usepackage{amssymb, latexsym, amsmath}
    \newcommand{\setleftmargin}[1]{
      \addtolength{\textwidth}{\oddsidemargin}
      \addtolength{\textwidth}{1in}
      \addtolength{\textwidth}{-#1}
      \setlength{\oddsidemargin}{-1in}
      \addtolength{\oddsidemargin}{#1}
      \setlength{\evensidemargin}{\oddsidemargin}
    }
    \newcommand{\setrightmargin}[1]{
      \setlength{\textwidth}{8.5in}
      \addtolength{\textwidth}{-\oddsidemargin}
      \addtolength{\textwidth}{-1in}
      \addtolength{\textwidth}{-#1}
    }
    \newcommand{\settopmargin}[1]{
      \addtolength{\textheight}{\topmargin}
      \addtolength{\textheight}{1in}
      \addtolength{\textheight}{\headheight}
      \addtolength{\textheight}{\headsep}
      \addtolength{\textheight}{-#1}
      \setlength{\topmargin}{-1in}
      \addtolength{\topmargin}{-\headheight}
      \addtolength{\topmargin}{-\headsep}
      \addtolength{\topmargin}{#1}
    }
    \newcommand{\setbottommargin}[1]{
      \setlength{\textheight}{11in}
      \addtolength{\textheight}{-\topmargin}
      \addtolength{\textheight}{-1in}
      \addtolength{\textheight}{-\footskip}
      \addtolength{\textheight}{-#1}
    }
    \newcommand{\setallmargins}[1]{
      \settopmargin{#1}
      \setbottommargin{#1}
      \setleftmargin{#1}
      \setrightmargin{#1}
    } \setallmargins{1.2in}
\title{Generalized Lindemann-Weierstrass and Gelfond-Schneider-Baker Theorems
}
\author{
Suk-Geun Hwang, \thanks{
Department of Mathematics Education, Kyungpook University, Taegu 41566, Rep. of Korea, Email:
sghwang@knu.ac.kr}
Choon Ho Lee, \thanks{
Department of Mathematics,
College of Natural Sciences, Hoseo University,
20 79th-Gil Hoseo-Ro Baebang-Eup,
Asan Choong-Nam,
336-795,
Korea, Email: chlee@hoseo.edu} and
Ki-Bong Nam \thanks{ Dept. of Math., Univ. of Wisconsin-Whitewater,
Whitewater, WI 53190, USA,  Email: namk@uww.edu: The fourth author was partially supported by Strategic Priorities Fund  of the Dept. of Math., UW-W. } 
Rachel M Chaphalkar, \thanks
{ Dept. of Math., Univ. of Wisconsin-Whitewater,
Whitewater, WI 53190, USA,  Email: chaphalr@uww.edu}}
\begin{document}
\maketitle
 \begin{abstract}
We generalize Lindemann-Weierstrass theorem and Gelfond -Schneider-Baker Theorem.
We find new transcendental numbers in this work.
There are several methods to find transcendental numbers in the work.
Recently transcendental numbers are applicable for cryptography (\cite{G}, \cite{K}, \cite{V}). Since we are able to make
many tables of random numbers, the new transcendental 
numbers will be applicable for encryption and decryption in this work (\cite{V}, \cite{Z}). 
  \end{abstract}

      \newtheorem{lemma}{Lemma}
    \newtheorem{prop}{Proposition}
    \newtheorem{thm}{Theorem}
    \newtheorem{coro}{Corollary}
    \newtheorem{definition}{Definition}
\newtheorem{exa}{Example}[section]
\newtheorem{note}{Note}[section]


\section{Preliminaries}

Let ${\mathbb Z}$ be the set of all integers and ${\mathbb N}$ $(\hbox{resp}.$ ${\Bbb Z}^+)$ be the
set of all non-negative $(\hbox{resp. positive})$ integers. Let ${\mathbb Q}$ be the set of all rational numbers and ${\Bbb R}$ be the set of all real numbers.
Throughout the paper ${\Bbb F}_A$ denotes the field of all algebraic numbers.
We can define a degree on the mapping ring
${\Bbb Q}[f(x)]$ as the polynomial ring ${\Bbb Q}[x]$ (see \cite{C}).
The function ring ${\Bbb Q}[f(x)^{\pm 1}]$ is the Laurent extension of ${\Bbb Q}[f(x)]$. Similarly, we define a function ring with several variables.
The equations have different algebraic structures (see \cite{C}). 
Let $A$ be an additive subgroup of ${\Bbb F}$.
Let ${\mathcal T}$ be the set of all the transcendental numbers and let us define $\sim$ on ${\mathcal T}$ with respect to $A$.
For given two transcendental numbers $\tau_1$ and $\tau_2$, we define $\sim$ as follows:\\
$\tau_1\sim \tau_2$ if there is an element $r$ of $A$ such that $\tau_1=\tau_2+r$. Then $\sim$ is an equivalence relation (\cite{L1} and \cite{ST}). For a transcendental number
$\tau$ and $A$, ${\overline \tau}_A$ denotes the equivalence class of $\tau$ with respect to $\sim$. Let $\overline {\mathcal T}$ be the set of all equivalence classes of ${\mathcal T}$ with respect to 
$\sim$.




\section{Main Results}

Let us start with the well-known Lindemann-Weierstrass theorem \cite{N}:
\begin{thm}\label{LW}
For $\alpha_1,\cdots, \alpha_m\in {\Bbb F}_A$, if $\alpha_1,\cdots, \alpha_m$ are distinct non-zero algebraic numbers,
then $e^{\alpha_1}, \cdots, e^{\alpha_m}$ are linearly independent over ${\Bbb F}_A$ where at least one of
$\alpha_1,\cdots, \alpha_m$ is not zero.
\end{thm}

\noindent
Theorem \ref{LW} shows that for $c_1,\cdots, c_m\in {\Bbb Q},$
$c_1 e^{\alpha_1}+\cdots +c_me^{\alpha_m}$ is a transcendental number where
at least one of  $c_1,\cdots, c_m$ is not zero. We call the transcendental number $l$ as 
a Lindemann-Weierstrass number and we define ${\mathcal LW}$ as the set of all Lindemann-Weierstrass transcendental numbers of Theorem 1.
For $\tau_1,\tau_2\in {\mathcal T},$ both of $\tau_1+\tau_2$ and $\tau_1-\tau_2$ are generally not transcendental numbers.
But for any $\tau_1,\tau_2\in {\mathcal LW}$, $\tau_1\pm \tau_2\in {\mathcal LW}.$
By Theorem 1, every element of ${\Bbb Q}[e^{\alpha_1}, \cdots ,e^{\alpha_n}]$ is a transcendental number.
The following is Corollary 2.3 of the book \cite{Pa}:\\
Let $P(z_0,\cdots, z_s)\in {\Bbb F}_A[z_0,\cdots, z_s]$ be a non-zero
polynomial that is not divisible by any non-constant polynomial in $z_0.$ 
If $\alpha_1,\cdots ,\alpha_s$ are algebraic number that are linearly independent
over ${\Bbb Q},$ and $\tau$ is a zero of the equation
$P(z_0,\cdots, z_s)=0,$ then $\tau$ is transcendental.
The following results generalize Corollary 2.3 of the book \cite{Pa}.

\begin{thm}\label{GLW}
Let $f(x), f_1(x),\cdots , f_m(x), g_1(x),\cdots , g_m(x)\in {\Bbb F}_A[x]$ and let $\alpha_1,\cdots, \alpha_m$ be distinct non-zero algebraic numbers. If the equation
\begin{eqnarray}\label{eq1}
& &g_1(x) e^{\alpha_1f_1(x)}+\cdots + g_m(x) e^{\alpha_m f_m(x)}=f(x)
\end{eqnarray}
has a non-zero solution $\alpha$ such that there exists $i$ of $\{1,\cdots, m\}$ with
$f_i(\alpha)\neq 0\neq g_i(\alpha)$, and $f(\alpha)\neq 0$, then
the solution $\alpha$ is a transcendental number
where $f(x)$ is not zero. 
Every Lindemann-Weierstrass number is a solution of an equation (\ref{eq1}) of the form where $f(x)\in {\Bbb Q}[x]$.
\end{thm}
{\it Proof.}
Let us assume that $\alpha$ is a non-zero solution of the equation $g_1(\alpha)e^{\alpha_1f_1(x)}+\cdots +g_m(\alpha)e^{\alpha_m f_m(x)}=f(x)$ which 
is called a zero of $f(x)$.
If $\alpha$ is a transcendental number, then there is nothing to prove.
Let us assume that $\alpha$ is not a transcendental number. Then by Theorem 1 and Corollary 2.3 of \cite{Pa},
$g_1(\alpha)e^{\alpha_1f_1(\alpha)}+\cdots +g_m(\alpha)e^{\alpha_m f_m(\alpha)}$ is a transcendental number. But
$f(\alpha)$ is an algebraic number. So we have a contradiction. Therefore the solution $\alpha$ of the
equation is a transcendental number.
In the equation (\ref{eq1}), especially if $c_1e^{\alpha_1 x}+\cdots +c_me^{\alpha_m x}=f(x)$ and $f(x)=c_1e^{\alpha_1 }+\cdots +c_me^{\alpha_m },$
then $x=1$ which is where $c_1,\cdots ,c_m\in {\Bbb Q}$.
So every Lindemann-Weierstrass number is a solution of the equation (\ref{eq1}).
\quad $\Box$\\

\bigskip

\noindent 
Let ${\mathcal {LWE}}$ be the set of all transcendental numbers which are in Theorem 2. Then the cardinality $|{\mathcal {LWE}}|$ of ${\mathcal {LWE}}$
is $\aleph_0$.
For any non-zero algebraic number $\alpha$, $\ln \alpha, \sin \alpha, \cos \alpha, \tan \alpha, \csc \alpha, $ $
\sec \alpha, \cot \alpha, $ $\sinh \alpha, $ $\cosh \alpha,$ $ \tanh \alpha,$ and $\coth \alpha$ are transcendental numbers (see Theorem 9.11 of \cite{N}). In addition
$\ln \alpha$ is transcendental number where $\alpha$ is a non-zero algebraic number (see Theorem 9.11 of \cite{N}).
\begin{coro}
Let $f(x)$ be one of $ \ln x, \sin x, \cos x, \tan x, \csc x, \sec x, \cot x, \sinh x,$ $ \cosh x,$ $ \tanh x,$ and $\coth x$. For $g\in {\Bbb Q}[f(x)]$ and $h(x)\in {\Bbb Q}[x]$,
if an equation $g=h(x)$ has a non-zero solution $\alpha$ such that $f(\alpha)\neq 0 \neq g(f(\alpha))$, then the solution $\alpha$ is a transcendental number.
\end{coro}
{\it Proof.}
Let $\alpha$ be a solution of the equation $g=h(x)$. If $\alpha$ is algebraic, then $h(\alpha)$ and $g(f(\alpha))$ algebraic. 
So the proof of the corollary is straightforward by Theorem 2. \quad $\Box$

\begin{coro}
Let $f(x)$ be one of $\ln x, \sin x, \cos x, \tan x, \csc x, \sec x, $ $\cot x, $ $\sinh x,$ $ \cosh x,$ $ \tanh x,$ and $\coth x$.
For a given equation $h_1(x)f(g(x))=h_2(x)$, if the equation has a non-zero solution $\alpha$ such that all of $h_1(\alpha), g(\alpha), f(g(\alpha)),$ and $h_2(\alpha)$ are not zeros, then the solution $\alpha$ is a transcendental number where
$g(x),h_1(x), h_2(x)\in {\Bbb Q}[x]$. 
\end{coro}
{\it Proof.}
The proof of the corollary is straightforward by the proofs of Theorem 2 and Corollary 1.
\quad $\Box$\\

\begin{coro}
Let $\hbox {arc} \sin x=f(x),$ $arc \cos x=f(x),$ $arc \tan x=f(x),$ $arc \cot x=f(x),$ $arc \sec x=f(x),$ $arc \csc x=f(x)$ be given well-defined equations
where $f(x)\in {\Bbb Q}[x].$ If one of those equations has a non-zero solution $\alpha$ such that $f(\alpha)\neq 0$, then the solution $\alpha$ is a transcendental number.
\end{coro}
{\it Proof.}
The proof of the corollary is straightforward by Theorem 2.
\quad $\Box$

\bigskip

\noindent {\bf Note 1.} \\
\noindent 
A real number solution of $e^x+x-12=0$ is a transcendental number, that is 
$2.27472787147\cdots.$ 
The equation $e^x+x-12=0$ has countably many complex number solutions. The equation $\pi^x+4x=49$ has a real number solution
which is the transcendental number $x=3.14097\cdots $. For solving the equation $e^x+x-12=0$, we need to determine an appropriate 
value of $e$, then solve the equation. A solution of the equation is highly dependent on the value of $e$ and a solution of the equation gives
a good random numbers. For the equation $e^x+x-12=0$, if we take an approximate value of $e$, then we get a very long 
rational number which is an approximation of a zero of $e^x+x-12.$ So this number gives a table of random numbers as well \cite{F}.
\quad $\Box$

\begin{coro}
Let $f(x)$ be one of $e^{x^j}, \sin x, \cos x, \tan x, \csc x, \sec x, \cot x, \sinh x, \cosh x,$ $ \tanh x,$ and $\coth x$ where $j$ is a fixed positive integer.
For a given equation $(f(x)+a_1)^k=g(x)$, if the equation has a non-zero solution $\alpha$ such that 
$f(\alpha)$ and $g(\alpha)$ are not zeros, then the solution $\alpha$ is a transcendental number where
$g(x)\in {\Bbb Q}[x]$, $a_1$ is a given algebraic number, and $k$ is a given positive integer. 
\end{coro}
{\it Proof.}
Let $\alpha$ be a real solution of the equation of the corollary. Since $\sqrt [k] {\alpha}$ is algebraic, the proof of the corollary is straightforward by the proof of Theorem 2. 
\quad $\Box$

\begin{coro}
Let $f(x)$ be one of $e^{x^j}, \sin x, \cos x, \tan x, \csc x, \sec x, \cot x, \sinh x, \cosh x, \tanh x,$ and $\coth x$. 
For any $g\in {\Bbb Q}[f(x)]$ and for $h(x)\in {\Bbb Q}[x]$,
if $1\leq \deg (g)\leq 4$ 
and an equation $g=h(x)$ has a non-zero solution $\alpha$ such that all of $g(f(\alpha))$, $f(\alpha)$, and $h(\alpha)$ are not zeros, then the solution $\alpha$ is a transcendental number.
\end{coro}
{\it Proof.}
If $\tau$ is a non-zero zero of the equation
$g=h(x)$, then by the solution of a linear, the quadratic, Cardano's, or Ferarri's formulas respectively and the proofs of Theorem 2, 
we are able to prove the corollary easily. Let us omit it.
\quad $\Box$

\bigskip

\noindent {\bf Note 2.} 
By Corollary 1, we know that for algebraic number $\sinh \alpha$ and $\cosh \alpha$ are transcendental numbers as well.
A real number solution of $xe^x=-x+12$ is between $1.7$ and $1.8.$ The equation $e^x-x+7=0$ has no real solution, but it has countably many complex transcendental number solutions. One of the complex number solutions of
the equation $e^x-x+7=0$ is $x= 1.7701\cdots +  2.669613 \cdots \it i $ and its complex conjugate $x= 1.7701\cdots - 2.669613 \cdots \it i $ 
is a transcendental solution of the equation which has the minimal modulus of all the complex solutions of the equation. The complex solutions
of $e^x-x+7=0$ are discrete. For some values of $a_1,a_2, j,$ and $k$, the equation
$(e^{a_1x^j} +a_2)^k=f(x)$ has countably many complex number solutions.
\quad $\Box$\\

\bigskip

\noindent For two given transcendental numbers $\tau_1$ and $\tau_2,$ it is easy to prove that
at least one of $\tau_1+\tau_2$ and $\tau_1-\tau_2$ is a transcendental number. In addition we have the following results.
\bigskip

\begin{prop}
For any given transcendental numbers $\tau_1$ and $\tau_2$, $\tau_1+\tau_2i$ and $\tau_1-\tau_2i$ are transcendental numbers.
Specifically, $e+\pi i$ and $e-\pi i$ are transcendental numbers.
\end{prop}
{\it Proof.}
Let us assume that $\tau_1+\tau_2i$ is an algebraic number. This implies that $\tau_1-\tau_2i$ is algebraic.
Assume these algebraic numbers to be $\alpha_1$ and $\alpha_2$ respectively.
By adding them, we have that $2\tau_1=\alpha_1+\alpha_2.$ This implies that $2\tau_1$ is a transcendental number but
$\alpha_1+\alpha_2$ is an algebraic number. We have a contradiction and so this implies that
$\tau_1+\tau_2i$ is a transcendental number. 
Symmetrically $\tau_1-\tau_2i$ is a transcendental number.
This completes the proof of the proposition.
\quad $\Box$\\


\begin{prop}
Let $A_1$ and $A_2$ be additive subgroups of the field ${\Bbb F}$ of all algebraic numbers. If $A_1\subset A_2,$
then for any $\tau \in {\mathcal T}$, ${\overline \tau}_{A_1}\subset {\overline \tau}_{A_2}$.
\end{prop}
{\it Proof.}
For any $\tau_1 \in {\overline \tau}_{A_1},$ $\tau_1=\tau+a_1$ for $a_1\in A_1\subset A_2.$ This implies that
$\tau_1\in {\overline \tau}_{A_2}.$ Thus we have that 
${\overline \tau}_{A_1}\subset {\overline \tau}_{A_2}$.
\quad $\Box$

\bigskip

The well-known Gelfond-Schneider-Baker theorem is the following \cite{B}:

\begin{thm}\label{GSB}
For $\alpha_1,\cdots, \alpha_n, \beta_1,\cdots, \beta_n\in {\Bbb F}_A$ with $\alpha_1, \cdots, \alpha_n\neq 0,1,$ if $\beta_1, \cdots, \beta_n$ is irrational, then $\alpha_1^{\beta_2}\cdots \alpha_n^{\beta}$
is a transcendental number.
\end{thm}

Let us call the transcendental number of Theorem \ref{GSB} a Gelfond-Schneider 
(transcendental) number and let ${\mathcal {GSB}}$ be the set of all the
Gelfond-Schneider-Baker numbers. Note that $|{\mathcal {GSB}}|=\aleph_1$. We have the following generalized
Gelfond-Schneider-Baker theorem.

\begin{thm}\label{GGSB}
For $f_1(x),\cdots, f_n(x), g(x)\in {\Bbb F}_A[x]$ and $\alpha_1,\cdots, \alpha_n, \beta_1, \cdots ,\beta_n\in {\Bbb F}_A$ with $\alpha_1,\cdots, $ $\alpha_n \neq 0,1,$ if $\beta_1, \cdots ,\beta_n$ are irrational and
\begin{eqnarray}\label{GS1}
& &(\alpha_1 f_1(x))^{\beta_1} \cdots (\alpha_n f_n(x))^{\beta_n} =g(x)
\end{eqnarray}
has a  non-zero solution $\alpha,$ then
$\alpha$ is a transcendental number where all of $f_1(\alpha),\cdots ,f_n(\alpha)$ are not zeros and $g(\alpha)$ is not zero.
\end{thm}
{\it Proof.}
The proof of the theorem is very similar to the proof of Theorem 2, so let us omit it.
\quad $\Box$\\

\bigskip

\noindent 
Let ${\mathcal {GSBE}}$ be the set of all transcendental numbers which are in Theorem 4. Then $|{\mathcal {GSBE}}|$ is $\aleph_1$.
Especially, if we put $f_1(x)=\cdots =f_n(x)=1$ and $g(x)= \alpha_1^{\beta_1}\cdots \alpha_n^{\beta_n}$, then
every element of ${\mathcal {GSB}}$ is a transcendental number of Theorem \ref{GGSB}.
Solutions of the equation $(3x)^{\sqrt 7}=x^2+10x+5$ are $x=0.932103\cdots$ and $-0.395261\cdots -0.173148\cdots i$
which are transcendental numbers.

\bigskip

\noindent {\bf Further Consideration:}\\

\noindent 
For given transcendental numbers $\tau_i$, $1\leq i\leq 4,$ is the number 
$\tau_1+\tau_1 {\mathbf i}+ \tau_2 {\mathbf j}+\tau_3 {\mathbf k}$ a transcendental number in the 
quaternions ${\mathbf Q}$ (see \cite{E})?
Since transcendental numbers such as $e$ can be used to make a key for encryption and decryption [11],  
we may use additional transcendental numbers found through this work for the same purpose, strengthening the key \cite{K}.


\bigskip



\begin{thebibliography}
\frenchspacing
\parskip 11pt


\bibitem{B} A. Baker, {\it Transcendental Number
Theory}, Springer-Verlag, 1982.

\bibitem{C} Rachel M Chaphalkar, 
Suk-Geun Hwang, 
Choon Ho Lee, and
Ki-Bong Nam {\it New Gelfond-Type Transcendental Numbers}, Math arXiv:2106.04055v2, 2021.\\


\bibitem{E} Samuel Eilenberg and Ivan Niven, {\it The fundamental theorem of algebra for quaternions}, 
Bull. Amer. Math. Soc. 50(4): 246-248 (April 1944).


\bibitem{F} Pierre-Alain Fouque, J. Stern, J. Wackers, {\it CryptoComputing with Rationals}, Financial Cryptography, 2002.
Computer Science


\bibitem{G} V. M. Silva Garcia, M. D. Gonzalez, R. F. Carapia, E. Vega-A,
{\it A Novel Method for Image Encryption Based on Chaos and Transcendental Numbers},
IEEE, Vol. 7, 2019.



\bibitem{K} R. M. Kumar, S. S. P. Kumar, {\it Secure Encryption/Decryption Technique using Transcendental Number},
International Journal of Computer, 2015.



\bibitem{L1} William J. LeVeque, {\it Fundamentals of Number Theory},
Cambridge University Press, London, 1975.

\bibitem{N} Ivan Niven, {\it Irrational Numbers}, The Carus Mathematical Monographs,
Number 12, 1967.

\bibitem{Pa} A. N. Parshin and I. R. Shafarevich,  {\it Number Theory IV, Transcendental Numbers},
Springer Verlagl, 1998.


\bibitem{Sh} A. B. Shidlovskii, {\it Transcendental Numbers}, Translated from Russian by N. Koblitz, Walter de Gruyter, Berlin New York, 1989.


\bibitem{ST} Ian Stewart, {\it Algebraic Number Theory},
Chapman and Hall, 1979.

\bibitem{V} M. K. Viswanath, {\it Numbers and Cryptography},
www.m-hikari.com, 2014.

\bibitem{Z} Jiazhen Zhou,
Suk-Geun Hwang, 
Ki-Suk Lee,
Ki-Bong Nam, {\it Random Insertion Method based Security Protocols Using Transcendental Number Generators}, Preprint, 2022.


\end{thebibliography}
\end{document}